\documentclass{amsart}
\usepackage{amssymb,
amsmath,
amsthm,
graphicx,
multind,
eufrak,
xypic
}
\theoremstyle{plain}
\newtheorem{thm}{Theorem}[section]
\theoremstyle{definition}
\newtheorem{de}[thm]{Definition}
\newtheorem{ex}[thm]{Example}
\newtheorem{re}[thm]{Remark}

\newcommand{\CC}{{\mathbb C}}
\newcommand{\RR}{{\mathbb R}}
\newcommand{\NN}{{\mathbb N}}
\newcommand{\PP}{{\mathbb P}}
\renewcommand{\AA}{{\mathbb A}}

\newcommand{\pafg}[2][]{\frac{\partial #1}{\partial #2}}
{\begin{figure} \begin{center}}%
{\end{center} \end{figure}}

\newcommand{\codim}{\operatorname{codim}}
\newcommand{\ord}{\operatorname{ord}}

\newcommand{\bx}{\mathbf{x}}
\newcommand{\la}{\langle}
\newcommand{\ra}{\rangle}
\newcommand{\Der}{\operatorname{Der}\nolimits}

\newcommand{\End}{\operatorname{End}\nolimits}
\newcommand{\Hom}{\operatorname{Hom}\nolimits}
\newcommand{\U}{\operatorname{U}\nolimits}
\newcommand{\ad}{\operatorname{ad}\nolimits}
\newcommand{\liea}[1]{\mathfrak{#1}}
\newcommand{\lieg}[1]{\mathrm{#1}}

\begin{document}

\title{Transitive Lie algebras of vector fields---an overview}
\author[J.~Draisma]{Jan Draisma}
\address[Jan Draisma]{
Department of Mathematics and Computer Science\\
Technische Universiteit Eindhoven\\
P.O. Box 513, 5600 MB Eindhoven, The Netherlands\\
and Centrum voor Wiskunde en Informatica, Amsterdam,
The Netherlands}
\thanks{The author is supported by a Vidi grant from
the Netherlands Organisation for Scientific Research (NWO).}
\email{j.draisma@tue.nl}
\maketitle

\begin{abstract}
This overview paper is intended as a quick introduction to Lie algebras
of vector fields. Originally introduced in the late 19th century by
Sophus Lie to capture symmetries of ordinary differential equations,
these algebras, or {\em infinitesimal groups}, are a recurring theme in
20th-century research on Lie algebras. I will focus on so-called {\em
transitive} or even {\em primitive} Lie algebras, and explain their
theory due to Lie, Morozov, Dynkin, Guillemin, Sternberg, Blattner, and
others. This paper gives just one, subjective overview of the subject,
without trying to be exhaustive.
\end{abstract}

\maketitle

\section{Formal power series, vector fields, and Lie algebras}

In this section we introduce the basic objects of study: Lie algebras
of formal vector fields. The section concludes with a reading guide to
the remaining sections.

Throughout this paper, $K$ will denote a field of characteristic zero,
which we will specialise to $\RR$ or $\CC$ in some convergence questions.
Fix a natural number $n$, and let $K[[\bx]]:=K[[x_1,\ldots,x_n]]$ denote
the $K$-algebra of formal power series in the variables $x_1,\ldots,x_n$,
i.e., series of the form $\sum_{\alpha \in \NN^n} c_\alpha \bx^\alpha$,
where $\bx^{\alpha}$ denotes the monomial $x_1^{\alpha(1)} \cdots x_n^{\alpha(n)}$ and where
there are no restrictions on the coefficients $c_\alpha \in K$.
Every such series is of the
form $f=\sum_{d=0}^\infty f_d$ with $f_d$ a homogeneous polynomial of
degree $d$ in the variables $x_1,\ldots,x_n$. The minimal $d$ for which
$f_d$ is nonzero is called the {\em order} $\ord f$ of $f$. Let $M$
denote the subspace of $K[[\bx]]$ of elements of (strictly) positive
order. If $\ord f>0$, then $f$ can be written (in a non-unique way) as
$\sum_{i=1}^n x_i g_i$ where the minimum of the orders of the $g_i$ equals
$\ord f-1$.  Iterating this argument (with the $g_i$
instead of $f$ in the next iteration) shows that if $\ord f > d \geq 0$,
then $f$ is of the form $\bx^{\alpha_1} h_1+\ldots+\bx^{\alpha_N} h_N$
with monomials $\bx^{\alpha_i},\ \alpha_i \in \NN^n$ of degree $d$ and
all $h_i$ of order at least $1$. In particular, $f$ is then an element
of $M^{d+1}$, the space of (finite) linear combinations of products
of $d+1$ elements of $M$. This shows that in the metric on $K[[\bx]]$
defined by $d(f,g)=2^{-\ord(f-g)}$ the closed ball around $0$ with
radius $2^{-(d+1)}$ is contained in $M^{d+1}$; the opposite inclusion
is even easier.

By a {\em formal vector field (in $n$ variables)} over $K$ we mean
a {\em derivation} $X$ of the (commutative) formal power series
ring $K[[x_1,\ldots,x_n]]=:K[[\bx]]$, that is, a $K$-linear map
$K[[\bx]] \to K[[\bx]]$ satisfying Leibniz's rule $X(fg)=X(f)g+fX(g)$.
Set $f_i:=X(x_i)$. Then it turns out that $X$ equals the derivation
$Y:=\sum_{i=1}^n f_i \pafg{x_i}$ defined in the natural way. Indeed,
by Leibniz's rule and linearity the derivation $Z:=X-Y$ vanishes on
polynomials.  Moreover, by Leibniz's rule, $Z$ maps $M^d$ into $M^{d-1}$,
hence $Z$ is continuous with respect to the metric above. Since the
polynomials are dense in $K[[\bx]]$, we find that $Z$ is identically
zero, hence $X$ equals $Y$ as claimed.

The {\em order} $\ord X$ of a derivation
$X$ is defined as $-1+\min_i \ord X(x_i)$. Thus we have $\ord X(f) \geq
\ord X+\ord f$, and the elegance of this relation explains the term $-1$
in the definition of $\ord X$. Note that the minimal possible order of a
derivation is $-1$, for derivations of the form $\sum_{i} c_i \pafg{x_i}
+ \text{ higher order terms}$ with not all $c_i \in K$ equal to $0$.
The $K$-space $\Der K[[\bx]]$ of all formal vector fields forms a Lie
algebra with respect to the bracket $[X,Y]:=X \circ Y-Y \circ X$, and we
have $\ord[X,Y] \geq \ord X+\ord Y$. In particular, the derivations of order
at least $0$ form a subalgebra of $\Der K[[\bx]]$. 

We now come to the
central definition of this paper.

\begin{de}
A {\em Lie algebra of vector fields in $n$ variables} is a subalgebra $L$
of $\Der K[[\bx]]$, and it is called {\em transitive} if the subspace of
$L$ consisting of elements of non-negative order in $L$ has codimension
$n$ in $L$.
\end{de}

\begin{re}
This is equivalent to the condition that for every vector 
$(c_1,\ldots,c_n) \in K^n$ there is a vector field in $L$ of
the form $\sum c_i \frac{\partial}{\partial x_i} + \text{ higher order terms}$.
\end{re}

By the above, this subspace $L_0$ is then a sub{\em algebra}. Note that
$n$ is the maximal possible codimension of $L_0$ in $L$: any $n+1$-tuple
of derivations has a suitable linear combination $X$ where all of the
coefficients of the lowest-order terms $\pafg{x_1},\ldots,\pafg{x_n}$
cancel out, so that $X$ has order at least $0$.

To justify the terminology {\em transitive} assume that $K=\RR$ and
that all elements $X$ of $L$ are {\em convergent} in the sense that
all $X(x_i)$ are convergent power series near the origin. Then each $X$ can be
integrated, and transitivity means that the union, over all $X \in L$,
of the integral curves through the origin in $\RR^n$ contains an open
neighbourhood of the origin. Hence the corresponding {\em infinitesimal
group} acts locally transitively near the origin, and $L_0$ is the Lie
algebra of the isotropy group of the origin.

Lie's work and the subsequent work by various other authors to be
described below has as a principal aim the {\em classification} of
Lie algebras of vector fields {\em up to coordinate changes}. In
our present set-up, this means the following. If $\phi$ is an
automorphism of the $K$-algebra $K[[\bx]]$ and $X$ is a derivation,
then $X^\phi:=\phi X \phi^{-1}$ is again a derivation. Moreover, we have
$[X,Y]^\phi=[X^\phi,Y^\phi]$, so that if $L$ is a Lie algebra of vector
fields, then so is $L^\phi$. Thus the classification problem becomes:
{\em classify subalgebras of $\Der K[[\bx]]$ up to automorphisms}.

\begin{ex} \label{ex:DimensionOne}
Let us do the classification of {\em finite-dimensional} Lie algebras of
vector fields in {\em one} variable $x$. Let $L$ be such a Lie algebra,
and let $S \subseteq \{-1,0,1,2,\ldots\}$ be the set of orders of elements
of $L$. Clearly $|S|=\dim L$, so $S$ is finite.

Now note that the commutator of $X=(x^d + \text{h.o.t.})\pafg{x}$
(of order $d-1$) and $Y=(x^d + \text{h.o.t})\pafg{y}$ (of order $e-1$)
equals $((e-d)x^{d+e-1}+\text{h.o.t.})$, and that this is of order
$(d-1)+(e-1)$ if $d$ and $e$ are not equal. Hence if $2 \leq d <
e$, then the vector fields $Y,[X,Y],[X,[X,Y]],\ldots$ have orders 
\[ e-1<(d-1)+(e-1)<2(d-1)+(e-1)<\ldots. \]
We conclude that $S$ contains at most one positive number. 

Let $X \in L$ be of minimal order $d$ in $L$. There are three cases to
be considered:
\begin{description}
\item[$d>0$] By the argument above, $S=\{d\}$ and $L = \la X \ra_K$,
and the classification boils down to the classification of a single
vector field. In one variable this is straightforward: after a coordinate
change, $L$ equals $\la x^d \pafg{x} \ra$ (and for distinct $d$s these
Lie algebras are not equivalent).

\item[$d=0$] After a coordinate change, we may assume that
$X=x\pafg{x}$. Let $\ad(X)$ denote the linear map $L \to L,\ Y \mapsto
[X,Y]$. If $Y=(\sum_{i=2}^\infty c_i x^i)\pafg{x}$ is a second element
of $L$, then $\ad(X)^k Y=\sum_{i=2}^\infty (i-1)^k c_i \pafg{x}$.
Using Vandermonde determinants one sees that if $c_i$ is non-zero for
infinitely many $i$, then these vector fields for $k=0,1,2,3,\ldots$
are linearly independent, which contradicts the assumption that $L$ is
finite-dimensional. Hence only finitely many of the $c_i$ are non-zero,
and then the same Vandermonde argument shows that $x^i \pafg{x} \in L$
for all $i$ with $c_i$ non-zero. By the argument in the previous case,
there can only be one such $c_i$. Hence $L = \la x\pafg{x},x^i \pafg{x}
\ra_K$ for a single $i>1$, and for distinct $i$ these Lie algebras are
not equivalent under coordinate changes (although they are isomorphic
as abstract Lie algebras).

\item[$d=-1$] In this case, for any $Y \in L$ of non-negative order, we
have $\ord [X,Y]=\ord(Y)-1$. Hence $S$ is an interval $\{-1,0,\ldots,e\}$
for some $e \leq 1$. If $e=-1$, then $L$ is spanned by $X$, and after
a coordinate change we have $L=\la \pafg{x} \ra_K$. If $e=0$, then $L$
contains an element $Y$ of order $0$.  After a coordinate change we may
assume that $Y=x\pafg{x}$, and then an argument as in the previous case
shows that $L=\la \pafg{x},x\pafg{x} \ra_K$. Finally, if $e=1$, then we
may again assume that $L$ contains $x\pafg{x}$, and by a similar argument
as before we have $L=\la \pafg{x},x\pafg{x},x^2\pafg{x} \ra_K$. 
\end{description}
Among the Lie algebras just found, only those with $d=-1$ are
transitive. The largest one among them is isomorphic to $\liea{sl}_2(K)$,
the algebra of trace-zero $2 \times 2$-matrices with the commutator
as Lie bracket. In fact, it can also be derived as follows
(see, for instance, \cite[Example 1.58]{Olver86}): let the
group $\lieg{SL}_2$ act on the projective line $\PP^1 K$ via M\"obius
transformations. This gives rise to a map from $\liea{sl}_2$ into
vector fields on the usual affine chart $\AA^1 K \subseteq \PP^1 K$, and
with a suitable choice of coordinate these vector fields are 
$E=-\pafg{x},H=-2x\pafg{x},F=x^2\pafg{x}$, satisfying the familiar
commutation rules
\[ [H,E]=2E,\ [H,F]=-2F,\ [E,F]=H. \]
Elements with these commutation rules are called a {\em Chevalley basis}
of $\liea{sl}_2$. The subalgebra spanned by $E$ and $H$ corresponds
to the action of the group of upper triangular matrices (scaling and
translation but no inversion, in M\"obius terminology), and the subalgebra
spanned by the $E$ corresponds to the group of translations only.
\end{ex}

In a larger number of variables, the problem of classifying {\em all}
Lie algebras $L$ of vector fields up to coordinate changes turns out to
be intractable, already for the simple reason that classifying a {\em
single} vector field vanishing at the origin up to
coordinate changes (without, for instance,
any restrictions on its linear part as in Poincar\'e-Dulac theory) 
is intractable. 
Hence we will from now on restrict
our attention to {\em transitive} algebras $L$, whose structure turns
out to be very beautiful. The remainder of this paper is organised as
follows. In Section~\ref{sec:Lie} we describe Lie's classification
of transitive Lie algebras of vector fields in two variables, and
their appearance as symmetry algebras of scalar ordinary differential
equations. In Section~\ref{sec:GSB} we present the beautiful work by
Guillemin-Sternberg and by Blattner that relates transitive Lie algebras
to pairs of a Lie algebra and a subalgebra. In Section~\ref{sec:MD} we
describe the classification of the {\em primitive} ones among these pairs,
largely due to Morozov and Dynkin.  Finally, in Section~\ref{sec:Nice}
we go back to an observation by Lie on Lie algebras of vector fields with
``nice'' coefficients, which inspires two challenging
research problems that conclude the paper.

\section{Lie: transitive Lie algebras in dimension two}
\label{sec:Lie}

Towards the end of the 19th century, Lie classified finite-dimensional
transitive Lie algebras of vector field in {\em two} variables
\cite{Lie24}.  To achieve this, he divided these Lie algebras into two
classes, namely, the {\em primitive} ones and the {\em imprimitive}
ones. In modern Lie-algebraic terminology these notions are defined
as follows.

\begin{de}
Let $L$ be a transitive subalgebra of $\Der K[[x_1,\ldots,x_n]]$, and let
$L_0$ be its subalgebra of elements with non-negative order. Then $L$
is called {\em imprimitive} if there exists a subalgebra $L'$ of $L$
that lies strictly between $L_0$ and $L$.
\end{de}

As with transitivity, this terminology is best explained in the setting
where $K=\RR$ and the elements of $L$ are convergent. If $L'$ is a
subalgebra strictly between $L_0$ and $L$, of codimension $m,\ 0<m<n$
in $L$ say, then integrating $L'$ (using a weak version of Frobenius's
theorem on vector fields in involution; see for instance \cite[Theorem
1.40]{Olver86}) yields a foliation of a neighbourhood of $0$ in $\RR^n$
with sheets of dimension $n-m$, and these sheets are permuted by all
(local) one-parameter subgroups corresponding to elements in $L$. Thus
these sheets form a system of {\em imprimitivity}.
In the case of two variables, a subalgebra $L'$ as above corresponds to a
foliation of a neighbourhood of the origin in $\RR^2$ by curves permuted
by the infinitesimal group corresponding to $L$ ({\em eine invariante
Kurvenschare} in Lie's terminology), and $L$ is primitive if no such
foliation exists. 

We now give Lie's classification of primitive and imprimitive Lie
algebras in two variables $x,y$, where we adopt Lie's notation
$p:=\pafg{x},q:=\pafg{y}$.

\begin{thm}[Lie, \cite{Lie24}]
Every finite-dimensional, primitive Lie algebra in two variables over
an algebraically closed field $K$ of characteristic zero equals one of
the Lie algebras in Table~\ref{tab:Primitive}, up to a formal coordinate
change.
\end{thm}

In fact, if $K=\CC$, then ``formal'' in this theorem can be replaced
by ``convergent''; see Remark~\ref{re:Convergence} below.

\begin{table}
\begin{center}
\begin{tabular}{rll}
\hline
Type&	Lie algebra & Label\\
\hline
(5)&	$\langle p, q, xq, xp-yq, yp\rangle$&	A3\\
(6)&	$\langle p, q, xq, xp-yq, yp, xp+yq\rangle$&	A2\\
(8)&	$\langle p, q, xq, xp-yq, yp, xp+yq, x^2p+xyq, xyp+y^2q
	\rangle$&	A1\\
\hline
&&\\
\end{tabular} 
\caption{The primitive Lie algebras in two variables.  The third column
of table \ref{tab:Primitive} contains the the label given to these
algebras in \cite{GonzalezLopez92} and \cite{Lie24}. The {\em type}
will be explained in the next section.}
\label{tab:Primitive}
\end{center}
\end{table}

\begin{thm}[Lie, \cite{Lie24}]
Every imprimitive (but transitive) Lie algebra in two variables over an
algebraically closed field $K$ of characteristic zero can be moved by a
formal coordinate change into one of the 16 (families of) Lie algebras
in Table~\ref{tab:NonPrimitive}.
\end{thm}

The classification over the real numbers is slightly more
involved; see \cite{GonzalezLopez92b}.

\begin{table}
\begin{center}
\begin{tabular}{rrll}
\hline
Type&	Case&	Realisation& Label\\
\hline
(1,1)&& $\langle p, x^i\exp(\alpha x)q\rangle$,&	B$\beta$1,D1,D2\\
	&&where $i=0, \ldots, r_{\alpha} \geq 0$&\\ 
	&&and $\alpha$ in a non-empty finite set&\\
(1,2)&& $\langle p, yq, x^i\exp(\alpha x)q\rangle$,&	B$\beta$2,C2\\
	&&where $i=0, \ldots, r_{\alpha} \geq 0$&\\ 
	&&and $\alpha$ in a non-empty finite set&\\
(1,3)&& $\langle p, q, 2yq, -y^2q\rangle$&	C5\\
\hline 
(2,1)&1&$\langle p, xp+q \rangle$\\
	&&For a different $\liea{h}_1$, this is (1,1).\\
(2,1)&2&$\langle p,xp,x^iq \rangle, \text{ where }i=0, \ldots, r
	\geq 0$&	B$\gamma$1\\
	&&For $r=0$ and a different $\liea{h}_1$, this is (1,2).&\\
(2,2)&1&$\langle p, xp-\lambda yq, x^iq \rangle, \text{ where } i=0, 
	\ldots,r \geq 0 \text{ and } \lambda \neq 0$&B$\gamma$2,C8,D3\\
(2,2)&2&$\langle p, xp+((r+1)y+x^{r+1})q, x^iq \rangle$,&B$\gamma$3\\ 
	&&where $i=0, \ldots, r \geq 0$&\\
(2,2)&3&$\langle p, xp, yq, x^iq \rangle, \text{ where } i=0, \ldots, r
	\geq 0$&	B$\gamma$4,C3\\
(2,3)&& $\langle p, xp, q, 2yq, -y^2q \rangle$&	C6\\
\hline
(3,1)&1&$\langle p, 2xp+q, -x^2p-xq \rangle$&	B$\delta$1\\
(3,1)&2&$\langle p, 2xp, -x^2p, q \rangle$&\\
	&&For a different $\liea{h}_1$, this is (1,3).\\
(3,1)&3&$\langle p, 2xp-q, -x^2p+xq, q \rangle$&	B$\delta$2\\
(3,2)&1&$\langle p, 2xp-2yq, -x^2p+(1+2xy)q \rangle$&	C9\\
(3,2)&2&$\langle p, 2xp+ryq, -x^2p-rxyq, x^iq \rangle$,& B$\delta$3\\
	&&where $i=0, \ldots, r \geq 1$&\\
(3,2)&3&$\langle p, 2xp+ryq, -x^2p-rxyq, yq, x^iq \rangle$,& B$\delta$4\\ 
	&&where $i=0, \ldots, r \geq 0$&\\
	&& For $r=0$ and a different $\liea{h}_1$ this is (2,3).&\\
(3,3)&&	$\langle p, 2xp, -x^2p, q, 2yq, -y^2 q \rangle$&C7\\
\hline
&&&\\
\end{tabular}
\caption{The non-primitive transitive Lie algebras in two variables.
This table differs from Lie's table in that the origin
$(0,0)$ is always a regular point. Moreover, whenever an $\liea{sl}_2$
occurs, its Chevalley basis is contained in the basis given in our
table. The third column can be used for translation between this table
and Lie's table in \cite{Lie24}. The {\em type} will be explained in
the next section.}
\label{tab:NonPrimitive}
\end{center}
\end{table}

\begin{re} \label{re:Redundancy}
There is some redundancy in Table~\ref{tab:NonPrimitive}. For instance,
the Lie algebra $L=\la p,xp+q \ra$ of type $(2,1)$, case 1, can also be
brought into a Lie algebra in the family of type $(1,1)$. To see this,
consider the coordinate change $x=v\exp(-u)$ and $y=-u$ with inverse
$u=-y$ and $v=x\exp(-y)$. Under this coordinate change, we have
\begin{align*} 
p&=\pafg[u]{x} \pafg{u} + \pafg[v]{x} \pafg{v}=\exp(u) \pafg{v}\\
q&=\pafg[u]{y} \pafg{u} + \pafg[v]{y} \pafg{v}=-\pafg{u}-v \pafg{v}\\
xp+q&=-\pafg{u}
\end{align*}
so that the coordinate change moves the Lie algebra $L$ into the member
of the family in the first row with $\alpha$ running through the set
$\{1\}$ and $r_1=0$. A similar computation moves the Lie algebra of
type $(3,1)$, case 2, into a Lie algebra of type $(1,3)$. This
explains the lack of labels in the third row for these two Lie
algebras. This redundancy will be explained in the next section.
\end{re}

In the following section we will recast these classifications in
modern terminology, explaining the first column in both tables in the
process. But first let us quickly review {\em why} Lie was interested
in classifying Lie algebras of vector fields in two variables. This
discussion will not be needed in the following sections.

Let $A=K[[x,y]][y',y'',\ldots]$ be the ring of polynomials in the
(algebraically independent) variables $y'=y^{(1)},y''=y^{(2)},\ldots$
with coefficients that are formal power series in the variables
$x,y=y^{(0)}$. An element of $A$ represents the left-hand-side of a scalar
ordinary differential equation, with $x$ playing the role of independent
variable and $y,y',\ldots$ playing the role of the dependent variable
and its derivatives with respect to $x$. On $A$ one defines the {\em
total derivative} with respect to $x$ as the derivation
\[ D_x=\pafg{x}+\sum_{i=0}^\infty y^{(i+1)}\pafg{y^{(i)}}. \]
A vector field $X$ in two variables $x,y$ has a unique extension, still
denoted by $X$, to a $K$-linear derivation on $A$ for which $[X,D_x]$
equals $Q D_x$ for some $Q \in A$. This extension is called the {\em
prolongation} of $X$ (for this short characterisation of the prolongation
see \cite[Lemma 2.4]{Oudshoorn02}). Since the effect of $X$ on power
series in $x$ and $y$ is already prescribed, and since the prolongation
is to satisfy Leibniz's rule, the prolongation is determined by its value
on $y',y'',\ldots$.  Writing $f:=X(x)$ and $g:=X(y)$ and using short-hand
notation such as $f_{xy}$ for $\frac{\partial^2 f}{\partial x \partial
y}$, we find
\begin{align}
Q &= (QD_x)x = [X,D_x](x) = X(1)-D_x(X(x))=-f_x-f_y y'
\notag\\
X(y')=&X(D_x(y))=D_x(X(y))+[X,D_x](y) \notag\\
=&g_x+g_y y'+Q y'=g_x+(g_y-f_x)y'-f_y(y')^2 \notag\\
X(y'')=&X(D_x(y'))=D_x(X(y'))+[X,D_x](y') \notag\\
=&g_{xx}+(2g_{xy}-f_{xx})y'+(g_{yy}-f_{xy})(y')^2+(g_y-f_x)y'' \notag\\
&-f_{xy}(y')^2-f_{yy}(y')^3-2f_y y'y''+Q y'' \notag\\
=&g_{xx}+(2g_{xy}-f_{xx})y'+(g_{yy}-2f_{xy})(y')^2+(g_y-2
f_x)y'' \label{eq:Xy2}\\
&-f_{yy}(y')^3-3f_y y' y'', \notag
\end{align}
etc. The idea behind this definition is that the prolongation of $X$
describes the vector field induced by $X$ on higher jet spaces, where $y$
is considered the dependent variable and $x$ the independent variable:
if $K=\RR$ and $t \mapsto (x(t),y(t))$ is an integral curve of $X$
with $\frac{\mathrm{d}x}{\mathrm{d}t}(t_0) \neq 0$, then near $t_0$
we may write $t$ and $y$ as functions of $x$, repeatedly differentiate
$y$ with respect to $x$, and rewrite those derivatives as functions
of $t$. Then $t \mapsto (x(t),y(t),y_{x}(t),y_{xx}(t),\ldots)$ is the
unique integral curve of the prolongation of $X$ that projects to that
of $X$. 

Now we come to the central connection between Lie algebras of vector
fields and ordinary differential equations.

\begin{de}
Let $P$ be an element of $A$ of the form $y^{(m)}-Q$ with $m \geq 2$ a
natural number and $Q$ an element of $A$ involving none of the variables
$y^{(i)}$ with $i \geq m$, so that $P=0$ is an {\em explicit} o.d.e. of
order $m \geq 2$. Define $L(P)$ as the $K$-space of all vector fields
$X \in \Der K[[x,y]]$ whose prolongation satisfies $P|X(P)$, that is,
the polynomial $X(P)$ in $y',y'',\ldots$ is divisible by $P$ in the
ring $A$. The space $L(P)$ is called the {\em algebra of Lie
(point) symmetries} of the ordinary differential equation $P=0$.
\end{de}
For non-explicit $P$, the definition of $L(P)$ is slightly more
complicated; see \cite[Definition 2.5]{Oudshoorn02} in an algebraic
setting and the standard reference \cite[Chapter 2]{Olver86} in the
smooth setting (including the case of multiple independent variables);
in what follows we assume that $P$ is as in the definition. An easy
computation shows that $L(P)$ is, indeed, a Lie algebra. The rationale
for this definition can again be found in jet spaces: A local solution
$y=y(x)$ of $P=0$ gives rise to a curve in higher jet spaces, and applying
the infinitesimal transformation corresponding to (the prolongation of)
$X$ to such a ``jet curve'' should yield a curve on which $P$ vanishes,
as well, i.e., another solution curve.

For $P$ explicit of order at least two, one can show that $L(P)$ is
a finite-dimensional Lie algebra over over $K$ (with a uniform bound
on the dimension; see \cite[Exercise 2.27]{Olver86} and references
there). Having a classification of these Lie algebras at hand, one may
try to solve $P=0$ by making a change of coordinates that brings $L(P)$
into normal form, and analysing the solutions of equations $P$ for which
$L(P)$ is in normal form. This beautiful proposal by Lie involves steps
that are difficult from an algebraic and algorithmic perspective. In fact,
whether $L(P)$ is transitive can be determined algorithmically and if it
is transitive, the relevant entry in Table~\ref{tab:Primitive} can be
determined algorithmically, as well \cite{Draisma01,Reid91a,Reid91b},
but the coordinate change bringing the Lie algebra in normal form can
be arbitrarily complicated, and it is not clear when coordinate changes
involving only ``elementary functions'' exist. Nevertheless, special cases
of this idea are presently used in computer software for solving o.d.e.s.

\begin{re}
Note that $L(P)$ may very well {\em not} be transitive. In fact, apart
from proposing the solution strategy above, Lie also strongly advocated
the use of a single known vector field in $L(P)$ for trying to solve
$P$. Numerous examples of o.d.e.s of orders one, two, and
three where this works can be found in \cite{Lie1891}.
\end{re}

We conclude this section with two
examples, one well-known ``direct'' example and one ``inverse'' example.

\begin{ex}
First consider the o.d.e. $P=0$ with $P=y''$. Let $X=f \pafg{x}+g
\pafg{y}$ be a vector field, and denote its prolongation by $X$, as
well. For $X$ to be in $L(P)$, the expression \eqref{eq:Xy2} for $X(P)$
must be a multiple of $P$. This is equivalent to the system 
\begin{align}
g_{xx} &= 0 & 2g_{xy}-f_{xx}&=0 \label{eq:DetSys}\\
f_{yy} &= 0 & g_{yy}-2f_{xy}&=0 \notag
\end{align}
of linear p.d.e.s for $f$ and $g$. This is the so-called {\em determining
system} for $L(P)$. For general $P$, this determining system is a
typically over-determined system of linear p.d.e.s for the coefficients
$f,g$ whose solution space is the algebra of Lie symmetries.  Such
systems can be analysed with Buchberger-type algorithms \cite{Reid91a,
Schwarz98}. In fact, the relevant {\em Janet Bases} for systems of linear
p.d.e.s predate Buchberger's {\em Gr\"obner bases}. For a good recent
overview of their theory and implementations see \cite{Robertz07}.

In our particular case, the solution space to the determining system is
spanned by $(f,g)=(x^2,xy)$ and $(f,g)=(xy,y^2)$ together with all pairs
$(f,g)$ of (affine-)linear polynomials. We conclude that $L(P)$ is the
Lie algebra of type (8) in Table~\ref{tab:Primitive}.  This fact has the
following geometric interpretation: the solutions of the o.d.e. $y''=0$
are all (non-vertical) straight lines, and the Lie algebra of type
$(8)$ is the Lie algebra coming from the action of $\lieg{SL}_3$ on the
projective plane by means of projective linear transformations. These
transformations permute the collection of all straight lines. Note that
here the infinitesimal (or at least local) character of the algebra of Lie
symmetries becomes apparent: clearly, a projective linear transformation
may well transform a non-vertical line in the plane into a vertical
one, or even to the line at infinity. But given one solution (strictly
speaking again local), there is an open neighbourhood of the identity
in $\lieg{SL}_3$ that maps map the given solution to other solutions.
\end{ex}

\begin{ex}
Let us try to find an explicit o.d.e. $P=0$ of order $m \geq 2$ such
that $L(P)$ is the Lie algebra of type (5) in Table~\ref{tab:Primitive}.
For this, we compute the prolongations of some of the vector fields in
$L(P)$: First, the prolongation of $\pafg{x}$ is just $\pafg{x}$, and
$\pafg[P]{x}$ can only be a multiple of $P$ if it is zero. Hence $P=0$
is an {\em autonomous} o.d.e. Second, the prolongation of $\pafg{y}$ is
just $\pafg{y}$, and we find that $P$ does not depend
(explicitly) on $y$. Third, the
prolongation of $x \pafg{y}$ equals $x \pafg{y} + \pafg{y'}$. Since $P$
does not depend on $y$ and since the order of $P$ is at least two we find
that $P$ does not depend on $y'$ either. Hence the only {\em second-order}
explicit o.d.e. whose algebra of Lie symmetries contains the vector fields
so far is $y''=0$, and we already know its algebra of Lie symmetries
to be the entire algebra of type (8) in Table~\ref{tab:Primitive}. In
particular, if the algorithm of \cite{Draisma01} says that the algebra of Lie
symmetries of a given second-order o.d.e. is of type (8), then it can be
transformed into $y''=0$ by a coordinate change (note that we are being
slightly sloppy here, since formal coordinate changes in $x$ and $y$
do not necessarily map polynomial o.d.e.s to polynomial o.d.e.s).

Fourth, an easy induction shows that for $X=y \pafg{x}$ 
the expression $X(y^{(m)}),\ m \geq 2$ is a linear combination
of the monomials 
\[ y^{(1)}y^{(m)},y^{(2)}y^{(m-1)},\ldots,
y^{(\lfloor \frac{m+1}{2}\rfloor)} y^{(\lceil \frac{m+1}{2}
\lceil)},\] 
each with a non-zero coefficient. Applying $X$ to
$P=y^{(m)}-Q(y^{(2)},\ldots,y^{(m-1)})$ and replacing $y^{(m)}$ by $Q$ in
the result leaves, for instance, a non-zero multiple of $y^{(2)}y^{(m-1)}$
that cannot be cancelled by any of the other terms. Thus there exists
no explicit {\em polynomial} o.d.e. whose algebra of symmetries equals
the algebra (5) in Table~\ref{tab:Primitive}.  On the other hand,
if one enlarges the algebra $A$ to include more complicated o.d.e.s,
then there do exist examples with $L(P)$ equal to the Lie algebra (5),
such as 
\[ P:=y^{(4)}-\frac{5(y^{(3)})^2}{3y^{(2)}}-(y^{(2)})^{5/3}=0. \]
To check this, one needs only verify that the prolongation of $y
\pafg{x}$ is in $L(P)$ and that the prolongation of $x \pafg{x}$ is
{\em not}---otherwise $L(P)$ would include Lie algebra (6) of
Table~\ref{tab:Primitive}.  
\end{ex}

\section{Guillemin-Sternberg and Blattner: the realisation
theorem} \label{sec:GSB}

We return to transitive Lie algebras of vector fields in $n$ variables. In
their ground-breaking work \cite{Guillemin64} Guillemin and Sternberg
extracted the following fundamental properties of such a Lie
algebra $L \subseteq
\Der K[[\bx]]$ and its isotropy algebra $L_0 \subseteq L$ consisting
of all vector fields in $L$ of non-negative order: First, $L_0$ has
codimension $n$ in $L$---this was the defining condition for $L$ to be
transitive. Second, the only ideal of $L$ contained in $L_0$ is $\{0\}$;
here an {\em ideal} of $L$ is a subspace $I$ with the property that the
Lie bracket $[X,I]$ lies in $I$ for all $X$ in $L$. This second property
can be seen as follows. Suppose that an $L$-ideal $I$ contains a non-zero
element $Y_0=\sum_{i=1}^n f_i \pafg{x_i}$ of non-negative order $d$. Let
$f_i$ be of minimal order $d+1 \geq 1$, and let $x_j$ be a variable such
that the degree-$(d+1)$ part of $f$ really depends on $x_j$.  Since $L$
is transitive, it contains an element $X$ of the form $\pafg{x_j}+Z$
with $\ord Z \geq 0$. Then
\[ Y_1:=[X,Y_0]=[\pafg{x_j}+Z,Y_0]=[\pafg{x_j},Y_0]+[Z,Y_0] \]
Here the order of $[Z,Y_0]$ is at least $d$, while $[\pafg{x_j},Y_0]$ is
by construction non-zero with order $d-1$. Hence $Y_1$ has order $d-1$,
and it belongs again to the ideal $I$. Continuing in this
manner, one finds that $I$ contains elements of order $-1$
and is therefore not contained in $L_0$. 

The striking insight by Guillemin and Sternberg is that the
two conditions on $(L,L_0)$ suffice to determine a transitive Lie algebra
in $n$ variables. More precisely, in what follows we will
work with so-called {\em transitive pairs}.

\begin{de}
A {\em transitive pair} is a pair $(\liea{g},\liea{h})$ of a Lie
algebra $\liea{g}$ over $K$ with a subalgebra $\liea{h}$ of finite
codimension in $\liea{g}$; this is then called the {\em codimension of
the transitive pair}. The pair is called {\em effective} if $\liea{h}$
contains no non-zero $\liea{g}$-ideal.  The pair is called {\em primitive}
if $\liea{h}$ is a maximal (strict) subalgebra of $\liea{g}$.  A {\em
realisation} of a transitive pair $(\liea{g},\liea{h})$ is a Lie algebra
homomorphism $\phi: \liea{g} \to \Der K[[x_1,\ldots,x_n]]$ such that
$n=\codim_{\liea{g}} \liea{h}$ and such that $\liea{h}$ consists exactly
of those elements of $\liea{g}$ that are mapped by $\phi$ to vector
fields of non-negative order.
\end{de}

\begin{re}
Transitive pairs are the infinitesimal counterparts of homogeneous
spaces $G/H$ where $H$ is a closed subgroup of a Lie group $G$.  Such
spaces, and in particular the {\em symmetric spaces} among them, are
ubiquitous in differential geo\-metry; see \cite{Helgason62}. More general
manifolds-with-extra-structure, locally modelled on homogeneous
spaces and in some sense generalising Riemannian manifolds, are the
object of study in \cite{Sharpe97}. Again, transitive pairs
$(\liea{g},\liea{h})$ serve as indispensable infinitesimal models.
\end{re}

Loosely speaking, Guillemin and Sternberg proved that realisations of
transitive pairs exist and are unique up to coordinate
changes; here is the precise statement.

\begin{thm}[Realisation Theorem, Guillemin-Sternberg \cite{Guillemin64}]
Every transitive pair has a realisation. Moreover, if $\phi_1,\phi_2$
are realisations of the same transitive pair $(\liea{g},\liea{h})$,
then there exists a unique coordinate change $\psi:K[[\bx]] \to
K[[\bx]]$ such that, for all $X \in \liea{g}$, we have $\phi_2(X)=\psi
\phi_2(X) \psi^{-1}$. Furthermore, the {\em kernel} $\phi^{-1}(0)$
of any realisation $\phi$ of the transitive pair $(\liea{g},\liea{h})$
is equal to the largest $\liea{g}$-ideal contained in $\liea{h}$. In
particular, the realisation is an embedding if and only if
the pair is effective.
\end{thm}

In view of this realisation theorem, Lie's classification in
the previous section can be understood as a classification of all
finite-dimensional and effective transitive pairs $(\liea{g},\liea{h})$
of codimension $2$, up to natural isomorphisms between such pairs. This
modern view on Lie's classification helps in understanding
(and even re-doing) it. Assume that the pair is not primitive (for primitive Lie
algebras, see the following section). Then there exists a subalgebra
$\liea{l}$ such that $\liea{g} \supset \liea{l} \supset \liea{h}$,
where the inclusions are proper. In general, the transitive pairs
$(\liea{g},\liea{l})$ and $(\liea{l},\liea{h})$ will no longer be
effective. Denote the largest $\liea{g}$-ideal contained in $\liea{l}$
by $\liea{i}$ and the largest $\liea{l}$-ideal contained in $\liea{h}$
by $\liea{j}$. The quotients $(\liea{g}/\liea{i},\liea{l}/\liea{i})$
and $(\liea{l}/\liea{j},\liea{h}/\liea{j})$ are both
effective pairs of codimension $1$, and hence each is isomorphic to one of the
three transitive pairs found in Example~\ref{ex:DimensionOne}. Hence in
total, there are nine possibilities for the pair $(\dim \liea{g}/\liea{i},
\dim \liea{l}/\liea{j})$. This pair is called the {\em type} of the triple
$(\liea{g},\liea{l},\liea{h})$.  In our table \ref{tab:NonPrimitive},
the type is recorded in the first column but $\liea{l}$ is suppressed.
Sometimes, various intermediate subalgebras $\liea{l}$ can be chosen,
and the resulting triples may be of different types. As a consequence,
an effective pair $(\liea{g},\liea{h})$ may occur several times in the
classification.  When this happens, there is a reference to an earlier
entry of the list, to which it is isomorphic. This explains
the redundancy observed in Remark~\ref{re:Redundancy}.

\subsection{Blattner's proof of the realisation theorem}

Guillemin and Sternberg's original proof involves finding both the
required formal vector fields and the coordinate change relating
two realisations recursively, ``coefficient by coefficient''. A few
years later, Blattner came up with a beautiful coordinate-free proof
\cite{Blattner69}, where the uniqueness follows from the universal
property of certain representations. I will now spend a few paragraphs
explaining the proof, using basic Lie algebra theory as can be found,
for instance, in \cite{Jacobson62}. Making his proof explicit leads to
a {\em Realisation Formula} below. This formula will not be used in the
next section, but it will be used in the last section of this paper.

To sketch Blattner's proof, we fix a transitive pair $(\liea{g},\liea{h})$
and we recall the concept of {\em universal enveloping algebra}
$U(\liea{g})$ of $\liea{g}$. This is an associative algebra built from
$\liea{g}$ in the free-est possible way, subject to the condition that
the Lie bracket between elements in $\liea{g}$ coincides with their
commutator in $U(\liea{g})$. Formally, $U(\liea{g})$ is the quotient
of the {\em tensor algebra} $\bigoplus_{d=0}^\infty \liea{g}^{\otimes
d}$ generated by the vector space $\liea{g}$ by the (two-sided) ideal
generated by elements of the form
\[ X \otimes Y-Y \otimes X-[X,Y], \]
with $X$ and $Y$ running through $\liea{g}$. This expresses that, in
$U(\liea{g})$, the element $X \otimes Y-Y \otimes X \in \liea{g}^{\otimes
2}$ equals the Lie bracket $[X,Y] \in \liea{g}^{\otimes 1}$ of the
elements $X,Y$. The inclusion of $\liea{g}=\liea{g}^{\otimes 1}$ into the
tensor algebra gives a map $\liea{g} \to U(\liea{g})$, and this linear map
is, in fact, injective. For ease of exposition, we assume that $\liea{g}$
is finite-dimensional, although this is not needed for Blattner's proof
(and although the theory of transitive infinite-dimensional Lie algebras
is very rich; see \cite{Guillemin70} or, for the real $C^\infty$-case,
\cite{Singer65}). Let $X_1,\ldots,X_m$ be a basis of $\liea{g}$. Then
an easy induction shows that, modulo the defining relations above,
any {\em monomial}
\[ X_{i_1} \otimes \cdots \otimes X_{i_d} \] 
equals a linear combination of such monomials with the additional
property that the index sequence is weakly increasing. Such monomials
will be called {\em ordered}. Indeed, this linear combination will start
with the signed ordered monomial $\pm X_{i_{\pi(1)}} \otimes \cdots
\otimes X_{i_{\pi(d)}}$ (with $\pi$ a permutation rendering the indices
weakly increasing) and continue with monomials of lower degree. Thus the
ordered monomials span the universal enveloping algebra. The celebrated
Poincar\'e-Birkhoff-Witt theorem (see, for instance, \cite[Chapter V,
Theorem 3]{Jacobson62}) states that they are, in fact, a {\em basis} of
this algebra. By taking for $X_1,\ldots,X_{m-n}$ a basis of $\liea{h}$
and for $Y_1:=X_{m-n+1},\ldots,Y_n:=X_m$ a basis of a complementary
subspace, one finds that every element in $U(\liea{g})$ can be written
in a unique way as
\[ \sum_{\alpha \in \NN^n} u_\alpha Y_1^{\alpha_1} \cdots
Y_n^{\alpha_n} \]
with only finitely many terms non-zero and with coefficients $u_\alpha$
in the universal enveloping algebra $U(\liea{h}) \subseteq U(\liea{g})$
of $\liea{h}$.  In the case where $\liea{h}$ is zero-dimensional, this
is just a restatement of the Poincar\'e-Birkhoff-Witt theorem. Note that
we have left out the $\otimes$-signs here, and just use juxtaposition
to denote multiplication in $U(\liea{g})$. This notational choice is
especially important since other tensor products will soon
play a role.

We will get back to this explicit description of $U(\liea{g})$ as
a left-$U(\liea{h})$-module a bit later. For the moment, we will
not need it---and this is the beauty of Blatt\-ner's proof. The
universal enveloping algebra of $\liea{g}$ plays the same role in the
representation theory of $\liea{g}$ as the the group algebra $\CC G$ of
a finite group $G$ plays in representation theory of $G$. In particular,
a $\liea{g}$-module (or a representation of $\liea{g}$ as a Lie algebra
over $K$) gives rise to a $U(\liea{g})$-module (or a representation of
$U(\liea{g})$ as an associative algebra over $K$), and {\em
vice versa}.  Now Blattner introduces 
\[ A:=\Hom_{U(\liea{h})}(U(\liea{g}),K) \]
as a coordinate-free version of the algebra $K[[\bx]]$ of formal power
series; let me explain. First, $K$ is the trivial left-$\liea{h}$-module,
in which every element of $\liea{h}$ acts as zero. As a consequence,
all monomials in $U(\liea{h})$ of positive degree act as zero, as well,
while its unit element $1 \in \liea{h}^{\otimes 0}=K$ acts as unit element
(as one usually requires of associative-algebra representations).  Second,
$U(\liea{g})$ is a $U(\liea{h})$-module by left multiplication. Third,
$A$ is the space of all $K$-linear maps $a:U(\liea{g}) \to K$ such that
$a(vu)=va(u)$ for all $v \in U(\liea{h})$ and $u \in U(\liea{g})$. This
is equivalent to the condition that $a(Yv)=0$ for all $Y \in \liea{h}$.

One makes $A$ into a commutative algebra
as follows. Let $\Delta:U(\liea{g}) \to U(\liea{g}) \otimes U(\liea{g})$
denote the unique algebra homomorphism determined by $\Delta(X)= X\otimes
1+1 \otimes X$ for $X \in \liea{g} \subseteq
U(\liea{g})$.\footnote{This is the
Lie algebra analogue of the natural homomorphism $\CC G \to \CC G \otimes \CC G$
extending $g \mapsto g \otimes g$ for finite groups $G$, which allows
one to see the tensor product of two representations of $G$ as a new
representation of $G$, as opposed to a representation of $G \times
G$.} Then we define $a \cdot b$ as the the composition of the following
sequence of maps:
\[ U(\liea{g}) \overset{\Delta}{\to} U(\liea{g})\otimes U(\liea{g})
\overset{a \otimes b}{\to} K \otimes K \to K, \]
where the last map is just the natural (multiplication) isomorphism $K
\otimes K \to K$. A straightforward computation shows that $(a \cdot
b)(Yv)=0$ for all $Y \in \liea{h}$, so that $a \cdot b \in A$. The
fact that $\cdot$ is commutative follows from the fact that for every
$u \in U(\liea{g})$, $\Delta(u) \in U(\liea{g}) \otimes U(\liea{g})$
is invariant under swapping the two tensor factors---the {\em
co-commutativity} of $\Delta$. Similarly, the fact that $\Delta$
is {\em co-associative}, i.e., satisfies $(\Delta \otimes 1)\circ
\Delta=(1 \otimes \Delta)\circ \Delta$, implies that $\cdot$ is
associative. 

The algebra $A$ will serve as our coordinate-free model for the algebra
of formal power series---we will describe an isomorphism below.  Hence we
want $\liea{g}$ to act by derivations on $A$; this goes as follows. Given
$X \in \liea{g}$ and $a \in A$, we define a new element $Xa \in A$
by $(Xa)(u)=a(uX)$.  Note that $(Xa)(vu)=a(vuX)=va(uX)=v(Xa)(u)$ for
$v \in U(\liea{h})$, so that $Xa$ does indeed lie in $A$. In this way,
$A$ becomes a $\liea{g}$-module, as a straightforward computation shows.
Moreover, we have Leibniz's identity $X(a \cdot b)=X(a) \cdot
b + a \cdot X(b)$ for all $a,b \in A$ and $X \in \liea{g}$.  Indeed,
fixing $u \in \U(\liea{g})$ and writing $\Delta u=\sum_i u_i \otimes w_i$,
we find $\Delta(uX)=\sum_i u_i X \otimes w_i+\sum_i u_i \otimes w_i X$,
and hence
\[ 
(X(a\cdot b))u=(a \cdot b)(uX)=\sum_ia(u_iX) b(u_i)+\sum_i a(u_i)b(u_iX)=
(Xa \cdot b + a \cdot Xb)u, 
\]
as claimed. 

Thus we have $\liea{g}$ acting by means of derivations on an algebra
$A$. Now identifying $A$ with formal power series gives the existence
of a realisation of the transitive pair $(\liea{g},\liea{h})$. This
identification goes as follows. Fix a basis $Y_1,\ldots,Y_n$ of a
vector space complement of $\liea{h}$ in $\liea{g}$. We have already
seen that every element of $U(\liea{g})$ has a unique representation as
$U(\liea{h})$-linear combination of ordered monomials $Y_1^{\alpha_1}
\cdots Y_n^{\alpha_n}$. An element of $A$ assigns to each such monomial
a number $c_\alpha(a)$, and the map
\[ A \to K[[x_1,\ldots,x_n]],\ a \mapsto \sum_{\alpha \in NN^n} 
	\frac{c_\alpha}{\alpha!} \bx^\alpha, \]
where $\alpha!:=\prod_i \alpha_i!$, is the required algebra isomorphism.
To see where the factor ($1$ divided by) $\alpha!$ comes from, let $a_\alpha$
denote the unique element of $A$ that assigns $1$ to $Y_1^{\alpha_1}
\cdots Y_n^{\alpha_n}$ and $0$ to all other ordered monomials in the
$Y_i$. Then a straightforward computation shows that
\[ (a_\alpha \cdot a_\beta)(Y_1^{\alpha_1+\beta_1} \cdots
Y_n^{\alpha_n+\beta_n}) = \binom{\alpha+\beta}{\alpha,\beta} \]
while the value on all other ordered monomials is zero. This implies
that 
\[ a_\alpha \cdot a_\beta = \binom{\alpha+\beta}{\alpha,\beta}
a_{\alpha+\beta}, \]
a set of identities (with varying $\alpha$ and $\beta$) shared by the
``divided monomials'' $\frac{\bx^\alpha}{\alpha!}$.  This explains
the factorials appearing in the denominator, but it does more than
that: in positive characteristic, {\em divided power series} are in
many ways more natural than ordinary power series (and their algebra
is no longer isomorphic to the algebra of formal power series).
Blattner's construction automatically yields this algebra of divided
power series.\footnote{It turns out that, in positive characteristic,
the Lie algebra of derivations of the algebra of divided power series
has finite-dimensional simple sub-algebras, and these {\em
Cartan-type} Lie algebras are the main source of additional simple Lie
algebras when passing from the classification in characteristic zero
to the classification in positive characteristic \cite[Chapters 2 and
4]{Strade04}.}

Going back to characteristic zero, this explicit isomorphism between $A$
and the algebra of formal power series gives the following realisation
formula.

\begin{thm}[Realisation Formula] \label{thm:RealisationFormula}
Let $(\liea{g},\liea{h})$ be a transitive pair. Let $Y_1,\ldots,Y_n$ be
a basis of a vector space complement of $\liea{h}$ in $\liea{g}$. Then
the map $\phi$ sending $X \in \liea{g}$ to
\[ \phi(X)=\sum_{i=1}^n \left(\sum_{\alpha \in \NN^n} 
	\frac{a_{e_i}(Y_1^{\alpha_1} \cdots Y_n^{\alpha_n} X)}{\alpha!} 
	\bx^\alpha \right)\pafg{x_i}, \]
where $e_i=(0,\ldots,0,1,0,\ldots,0)$ with the $1$ on the $i$th position,
is a realisation of $(\liea{g},\liea{h})$.
\end{thm}

So to compute the coefficient of $\bx^\alpha$ in the coefficient
of $\pafg{x_i}$ in $\phi(X)$, one expresses $Y_1^{\alpha_1} \cdots
Y_n^{\alpha_n} X$ as a
$U(\liea{h})$-linear combination of ordered monomials in the $Y_j$,
and divides the (constant part of the) coefficient of $Y_i=Y_1^0 \cdots
Y_{i-1}^0 Y_i^1 Y_{i+1}^0 \cdots Y_n^0$ by $\alpha!$. Note that $X
\in \liea{g}$ is mapped to a derivation of non-negative order if and
only if $a_{e_i}(X)$ is zero for all $i$, which in turn is equivalent
to the statement that $X$ lies in $\liea{h}$. So $\phi$ is, indeed,
a realisation of $(\liea{g},\liea{h})$. The fact that its kernel is the
largest ideal of $\liea{g}$ contained in $\liea{h}$ is straightforward.

\begin{ex}
Let $\liea{g}=\liea{sl}_2=\la F,H,E \ra$ with the usual commutation
relations, let $\liea{h}=\la F,H\ra$, and choose the complementary
basis vector $Y_1:=E$. The linear function $a_1$ takes $E$ to $1$ and
all other ordered monomials $F^i H^j E^k$ to $0$. 
To compute $\phi(E)$, we have to consider
powers $E^d E=E^{d+1}$, which under $a_1$ are mapped to $1$ if $d=0$
and to $0$ otherwise. Hence $\phi(E)=\pafg{x}$. To compute $\phi(H)$,
compute 
\[ E^d H=E^{d-1}HE+E^{d-1}[E,H]=E^{d-1}HE-2 E^d=\ldots=
HE^d - 2d E^d. \]
Under $a_1$ this maps to $0$ for $d=0$ or $d>1$ and to $-2$ for $d=1$.
Hence $\phi(H)=-2 x \pafg{x}$. Finally, compute 
\begin{align*} 
E^d F&=E^{d-1}FE+E^{d-1}H=E^{d-1}FE+HE^{d-1}-2(d-1)E^{d-1}\\
&=\ldots=FE^d + d HE^{d-1} - d(d-1)E^{d-1}; 
\end{align*}
under $a_1$ this maps to $0$ unless $d=2$, in which case it maps to
$-2$. Thus $\phi(F)=\frac{-2}{2!} x^2 \pafg{x}=-x^2 \pafg{x}$. Up to a
sign, this is exactly the realisation found in
Example~\ref{ex:DimensionOne}.
\end{ex}

So far we have sketched the proof of the existence of
realisations. The uniqueness part of the Realisation Theorem follows from a
universal property of the $U(\liea{g})$-module $A$. For this, first note
that $A$ comes equipped with a natural $U(\liea{h})$-module map $a \mapsto
a(1)$ into $K$ (which, in terms of formal power series,
corresponds to evaluation at zero). Now given any other $U(\liea{g})$-module $B$
together with a $U(\liea{h})$-module map $\beta$ into $K$, there is
a unique $U(\liea{g})$-module map $\alpha: B \to A$ that
makes the following diagram commute:
\[ 
\xymatrix{
	A \ar[r]^{a \mapsto a(1)} &	K\\
	B \ar@{-->}[u]^{\exists! \alpha\ [U(\liea{g})-\text{homomorphism}]}
		\ar[ur]_{\beta\ 
		[U(\liea{h})-\text{homomorphism}]} &
		} 
\]
Indeed, $\alpha$ is forced to map
$b$ to the element $u \mapsto \beta(ub)$ of $A$, and this map does
the trick.\footnote{This construction is dual to the more familiar induction of
representations from a smaller group to a larger group, and the fact that
$U(\liea{h})$-homomorphisms $B \to K$ are in one-to-one correspondence
with $U(\liea{g})$-homomorphisms $B \to A$ is the analogue of Frobenius
reciprocity in this setting.} 

Now a {\em second} realisation of $\liea{g}$ gives rise, through
the same fixed isomorphism of $A$ with the algebra of formal power
series, to a second homomorphism $\phi_2:\liea{g} \to \Der(A)$
in addition to the first homomorphism $\phi: \liea{g} \to \Der(A)$
that we started with.  Taking $B$ above equal to $A$ {\em as a space},
but with $U(\liea{g})$-action coming from $\phi_2$, and taking $\beta$
equal to the map sending $b \in B$ to $b(1)$, the above shows that there
is a unique $U(\liea{g})$-module map $\alpha: B \to A$ such that $(\alpha
b)(1)=b(1)$. This property implies, as above, that $\alpha b=u \mapsto
(ub)(1)$. Somewhat lengthy, but straightforward computations show
that $\alpha$ is an algebra isomorphism from $B=A$ to $A$ intertwining
the realisations $\phi$ and $\phi_2$ as required by the theorem. This
concludes our description of Blattner's {\em proof from the book}.

\begin{re} \label{re:Convergence}
If $K$ equals $\RR$ or $\CC$, then one can show that the Realisation
Formula yields a realisation consisting of vector fields with a positive
radius of convergence around the origin \cite{Draisma02}. Moreover, the unique
formal coordinate change mapping that realisation to any other convergent
realisation can be shown to be convergent, as well. Hence the
Realisation Theorem implies that the classification of transitive Lie
algebras of convergent vector fields in $n$ variables is the same as
that of effective transitive pairs $(\liea{g},\liea{h})$ where
$\liea{h}$ has codimension $n$ in $\liea{g}$.
\end{re}

\section{Morozov and Dynkin: primitive Lie algebras}
\label{sec:MD}

While transitive Lie algebras are classified, up to coordinate changes,
by effective transitive pairs, the classification of such pairs themselves
in codimensions larger than two is very elaborate. Lie claims
to have completed the case of codimension $3$, but did not bother to
publish the complete result (see \cite{Lie24})---although he usually
did not eschew lengthy computations in his books.  Beyond codimension
three, the classification of transitive Lie algebras may well remain
out of reach.

But, as we saw in last section's discussion of Lie's classification
in two variables, transitive pairs of codimension $m$ give rise to
(in general, non-unique) sequences of {\em primitive} pairs whose
codimensions add up to $m$. Thus it makes sense to try and classify at
least {\em these}. This classification turns out to be very beautiful.
Recall that a primitive pair is a transitive pair $(\liea{g},\liea{h})$
with $\liea{h}$ a maximal subalgebra in $\liea{g}$. Adding the requirement
that the pair be effective, i.e., that $\liea{h}$ does not contain any
non-zero $\liea{g}$-ideal, turns out to leave only few possibilities
for $\liea{g}$. Here is the first theorem in that direction.

\begin{thm}[Morozov \cite{Morozov39}] \label{thm:Morozov}
Suppose that $(\liea{g},\liea{h})$ is an effective primitive pair.
Then either $\liea{g}$ is simple, or else we are in one of the following
two situations:
\begin{enumerate}
\item $\liea{g}$ is the direct sum $\liea{k} \oplus \liea{k}$ of two
isomorphic simple Lie algebras and $\liea{h}$ is the diagonal
subalgebra consisting of all pairs $(X,X),\ X \in \liea{k}$; or
\label{it:DirSum}
\item $\liea{g}$ is the semi-direct product $\liea{h} \ltimes \liea{m}$
with $\liea{m}$ an irreducible and faithful $\liea{h}$-module equipped
with trivial Lie bracket: $[\liea{m},\liea{m}]=0$, and with $\liea{h}$
is semisimple plus a center of dimension at most $1$.
\label{it:SemiDir}
\end{enumerate}
Conversely, in the latter two cases $(\liea{g},\liea{h})$ is
primitive and effective.
\end{thm}

This theorem reduces the classification of effective primitive pairs
to that of maximal subalgebras $\liea{h}$ of {\em simple} Lie algebras
$\liea{g}$ (for which the pair $(\liea{g},\liea{h})$ is trivially
primitive). We will discuss that classification below, after sketching
a proof of Morozov's result.

First assume that $\liea{g}$ is not semisimple. A direct consequence
of this is that $\liea{g}$ has a non-zero Abelian ideal $\liea{m}$.
Then $\liea{m}$ is not contained in $\liea{h}$ by effectiveness, so
$\liea{h}+\liea{m}$ is a subalgebra of $\liea{g}$ strictly containing
$\liea{h}$. Hence by primitivity we have $\liea{g}=\liea{h}+\liea{m}$
as vector spaces. But then $\liea{h} \cap \liea{m}$, being closed
under taking brackets with $\liea{h}$ and under taking under (zero)
brackets with $\liea{m}$, is an ideal in $\liea{g}$, hence zero by
effectiveness. Hence $\liea{g}=\liea{h} \oplus \liea{m}=\liea{h}
\ltimes \liea{m}$ as claimed. Any subspace $\liea{i}$ of $\liea{m}$
closed under taking brackets with $\liea{h}$ leads to a subalgebra
$\liea{h}\ltimes \liea{i}$ of $\liea{g}$, hence $\liea{i}$ equals zero
or all of $\liea{m}$; this proves that $\liea{m}$ is an irreducible
$\liea{h}$-module. Finally, consider the kernel of the homomorphism
$\liea{h} \to \End(\liea{m})$ sending $X$ to $[X,.]$. This kernel is
a $\liea{g}$-ideal, hence zero. We conclude that the representation of
$\liea{h}$ on $\liea{m}$ is faithful, as required. The only Lie algebras
$\liea{h}$ with faithful, irreducible representations are those that
are semisimple with a center of dimension at most one (which then acts
by means of scalars on $\liea{m}$).

Next assume that $\liea{g}$ is semisimple, and write $\liea{g}$ as
the direct sum of a simple ideal $\liea{k}$ and a semisimple ideal
$\liea{l}$, with $[\liea{k},\liea{l}]=0$. Then $\liea{k},\liea{l}$
are both not contained in $\liea{h}$ by effectiveness, and hence
$\liea{g}=\liea{h}+\liea{k}=\liea{h}+\liea{l}$ by primitivity. Let
$\pi_{\liea{k}}$ denote the projection $\liea{g} \to \liea{k}$ with
kernel $\liea{l}$. Then the equality $\liea{g}=\liea{h}+\liea{l}$ shows
that the restriction of $\pi_{\liea{k}}$ to $\liea{h}$ is surjective
onto $\liea{k}$. On the other hand, the kernel of this restriction is
$\liea{h} \cap \liea{l}$, which is an ideal in $\liea{g}$ and hence zero
by effectiveness. Hence $\pi_{\liea{k}}$ restricts to an isomorphism
$\liea{h} \to \liea{k}$. For entirely the same reason, the projection
$\pi_{\liea{l}}$ onto $\liea{l}$ along $\liea{k}$ restricts to an
isomorphism $\liea{h} \to \liea{l}$. Hence $\liea{k}$ and $\liea{l}$
are isomorphic simple Lie algebras, and $\liea{h}$ sits inside their
direct sum as a diagonal subalgebra. This concludes the proof of the
first part of Morozov's theorem. The converse, that the pairs under
\eqref{it:DirSum} and \eqref{it:SemiDir} are effective and primitive,
is straightforward.

\begin{ex}
Take $\liea{h}=\liea{gl}_2$ and $\liea{m}$ equal to the
standard two-dimensional module. Then the pair $(\liea{h} \ltimes
\liea{m},\liea{h})$ corresponds to the primitive Lie algebra of Type
(6) in Table~\ref{tab:Primitive}.
\end{ex}

As mentioned before, Morozov's theorem reduces the classification of
(finite-dimensional) effective, primitive pairs to that of maximal
subalgebras of simple Lie algebras $\liea{g}$. The following theorem deals
with the case where the field is algebraically closed and $\liea{g}$ is
{\em classical}, that is, $\liea{g}$ is isomorphic to the special linear
algebra $\liea{sl}_n$, the special orthogonal algebra $\liea{so}_n$,
or the symplectic algebra $\liea{sp}_{2m}$. Let $V$ be the standard
representation of $\liea{g}$, that is, $V$ equals $K^n$ in the first
two cases and $V=K^{2m}$ in the last case. In the last two cases, $V$
is equipped with a bilinear form which is symmetric in the orthogonal
case and skew-symmetric in the symplectic case.  Now let $\liea{h}$
be a maximal subalgebra of $\liea{g}$. Then Dynkin classifies the
possibilities as follows \cite{Dynkin57a}:
\begin{description}
\item[$\liea{h}$ acts reducibly on $V$] in this case, one of the
        following holds.
        \begin{enumerate}
        \item $\liea{h}$ is a maximal parabolic subalgebra of $\liea{g}$.
                These are of the form
                \[ \liea{p}(V'):=\{g \in \liea{g} \mid gV' \subseteq
V'\}
                \]
                with $V'$ a proper subspace of $V$, totally
                isotropic in case $\liea{g}$ is orthogonal or
                symplectic. Moreover, if $\liea{g}=\liea{o}_{2m}$,
		then $\dim U \neq m-1$.
        \item $\liea{g}=\liea{o}(V)$, and $\liea{h}=\liea{o}(U)
                \oplus \liea{o}(U^\bot)$ for some non-degenerate $U$,
                $0 \subsetneq U \subsetneq V$.
        \item $\liea{g}=\liea{sp}(V)$, and
                $\liea{h}=\liea{sp}(U) \oplus \liea{sp}(U^\bot)$ for
                some non-degenerate (and hence even-dimensional) $U$,
                $0 \subsetneq U \subsetneq V$.
        \end{enumerate}
\item[$\liea{h}$ acts irreducibly on $V$] then there are two
        possibilities.
        \begin{description}
        \item[$\liea{h}$ is not simple] then one of the following
        holds.
        \begin{enumerate}
        \item $\liea{h}=\liea{sl}(V)$, where $V \cong V_1 \otimes V_2$,
        and $\liea{h}=\liea{sl}(V_1) \oplus \liea{sl}(V_2)$. Here
        $V_1, V_2$ have dimensions $\geq 2$.
        \item $\liea{g}=\liea{sp}(V)$, where $V \cong V_1 \otimes V_2$,
        and $\liea{h}=\liea{sp}(V_1) \oplus \liea{o}(V_2)$.  Here $\dim
        V_1 \geq 2, \dim V_2 \geq 3$ but either $\dim V_2 \neq 4$ or
        $(\dim V_1,\dim V_2)$ equals $(2,4)$.  Moreover, $V_1$ is
        equipped with a non-degenerate skew bilinear form, and $V_2$
        with a non-degenerate symmetric bilinear form, such that the
        skew form on $V$ is the product of the two.
        \item $\liea{g}=\liea{o}(V)$, where $V \cong V_1 \otimes V_2$,
        and $\liea{h}=\liea{o}(V_1) \oplus \liea{o}(V_2)$.  Here $\dim
        V_1,\dim V_2 \geq 3$ but $\neq 4$. Moreover, each $V_i$ is
        equipped with a non-degenerate symmetric bilinear form, such
        that the symmetric form on $V$ is the product of the two.
        \item $\liea{g}=\liea{o}(V)$, where $V \cong V_1 \otimes V_2$,
        and $\liea{h}=\liea{sp}(V_1) \oplus \liea{sp}(V_2)$.  Here $\dim
        V_1,\dim V_2 \geq 2$. Moreover, each $V_i$ is equipped with a
        non-degenerate skew bilinear form, such that the symmetric form
        on $V$ is the product of the two.
        \end{enumerate}
        \item[$\liea{h}$ is simple] then one of the following holds.
        \begin{enumerate}
        \item $\liea{g}=\liea{sl}(V)$, and $\liea{h}=\liea{o}(V)$
        for some non-degenerate symmetric bilinear form on $V$.
        \item $\liea{g}=\liea{sl}(V)$, and $\liea{h}=\liea{sp}(V)$
        for some non-degenerate skew bilinear form on $V$.
        \item $\liea{g}=\liea{sl}(V)$, and $\liea{h}$ leaves invariant
        no bilinear form on $V$.
        \item $\liea{g}=\liea{o}(V)$, and $\liea{h}$ leaves invariant
        the symmetric bilinear form on $V$ defining $\liea{g}$.
        \item $\liea{g}=\liea{sp}(V)$, and $\liea{h}$ leaves invariant
        the skew bilinear form on $V$ defining $\liea{g}$.
        \end{enumerate}
        \end{description}
\end{description}
Conversely, still assuming that $\liea{g}$ is classical simple, if
$(\liea{g}, \liea{h})$ appears in the above list, then it is, as a rule,
a primitive pair. Only in the last three cases, there are exceptions to
this rule, and they are listed in Table 1 of \cite{Dynkin57a}.

This classification of maximal subalgebras of simple Lie algebras has
been implemented in the programme \verb+Lie+ \cite{LiE,maxsub}.  Dynkin
has also classified the maximal subalgebras of the {\em exceptional}
simple Lie algebras; see \cite{Dynkin57b}. By an elegant result due to
Karpelevich \cite{Karpelevich51}, one knows in advance that these are
all parabolic or reductive.

\section{Realisations with nice coefficients}
\label{sec:Nice}

This overview paper started off with Lie's problem of classifying
transitive Lie algebras of vector fields in $n$ variables. Through
the work of Guillemin-Sternberg and Blattner we have seen that
this is equivalent to classifying effective transitive pairs of Lie
algebras. Among these, the primitive pairs can actually be classified, as
we have seen in the previous section. 

\begin{figure}
\begin{center}
\includegraphics[width=.8\textwidth]{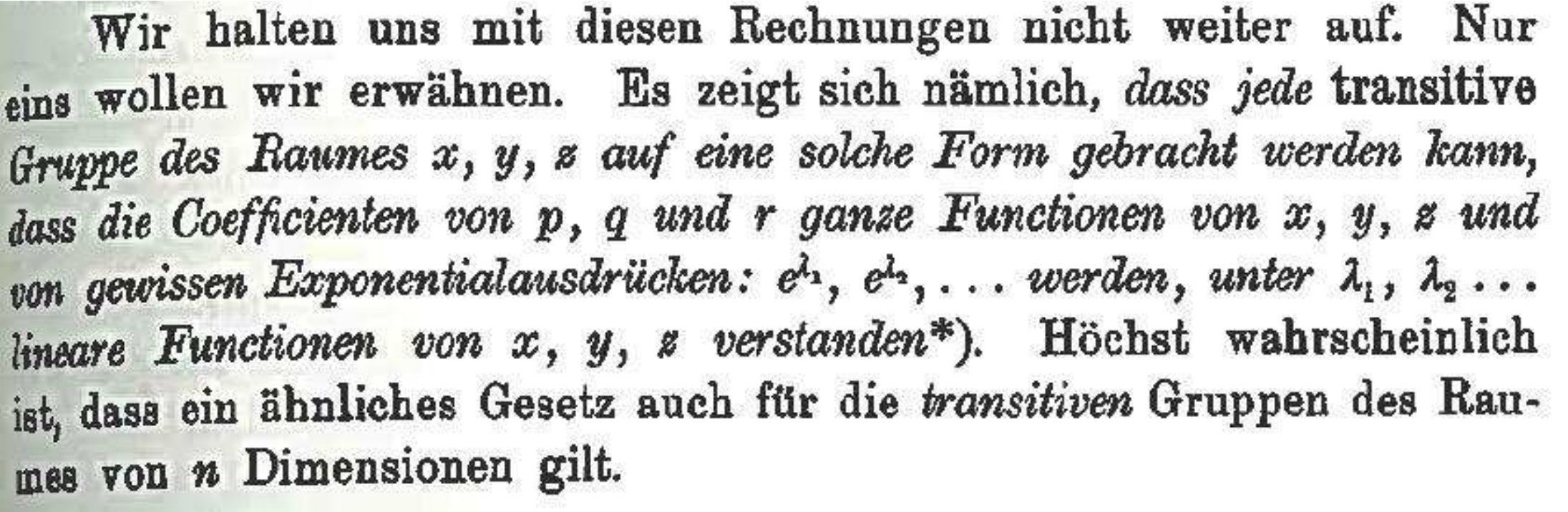}
\caption{Fragment from \cite{Lie1893}, page 177.}
\label{fig:exp}
\end{center}
\end{figure}

Now we will complete the circle as follows: given an effective transitive
pair $(\liea{g},\liea{h})$, we know that a realisation in terms of formal
power series exists, or even a realisation in terms of convergent power
series for $K=\RR$ or $\CC$---see Remark~\ref{re:Convergence}. But
sometimes we can actually find realisations with nicer coefficients,
such as polynomials, rational functions, or exponentials. Regarding such
realisations, Lie expressed the opinion in
Fragment~\ref{fig:exp}, translated as follows.

\begin{quote}
It turns out, {\em that every {\em transitive} group of $3$-space with
coordinates $x,y,z$ can be brought to a form in which the coefficients
of $p,q,$ and $r$ (Lie's notation for $\pafg{x}, \pafg{y},$ and
$\pafg{z}$---J.D.) are polynomial functions of $x,y,z,$ and certain 
exponential expressions $e^{\lambda_1}, e^{\lambda_2}, \ldots$, where
$\lambda_1, \lambda_2, \ldots$ denote linear functions of $x,y,z$.}
Very probably a similar statement holds for the {\em transitive} groups
of $n$-space.
\end{quote}

In fact, such realisations always exist when the pair is {\em very
imprimitive} in the sense of the following theorem.

\begin{thm}[Draisma, \cite{Draisma02}] \label{thm:VeryImprimitive}
Let $(\liea{g},\liea{h})$ be an effective transitive pair of codimension
$n$, and suppose that there exists a chain
\[ \liea{g}=\liea{g}_n\supseteq\liea{g}_{n-1} \supseteq \ldots \supseteq
\liea{g}_{1} \supseteq \liea{g}_0=\liea{h} \]
of subalgebras $\liea{g}_i$ with codimension $n-i$ for $i=0,\ldots,n$.
Then $(\liea{g},\liea{h})$ has a realisation $\phi: X \mapsto \sum_{i=1}^n
f_{i,X} \pafg{x_i}$ where for all $X$ and $i$ the coefficient $f_{i,X}$
is a polynomial in the $x_i$ and (finitely many) expressions $\exp(\lambda
x_i)$ with coefficients $\lambda$ in the algebraic closure of $K$.
\end{thm}

The proof of this theorem uses the Realisation Formula of
Theorem~\ref{thm:RealisationFormula}: as an ordered basis of a
complementary subspace to $\liea{h}$ one takes $Y_i \in \liea{g}_i
\setminus \liea{g}_{i-1}$ for $i=1,\ldots,n$, and then the coefficients
$f_{i,X}$ can be shown to satisfy non-trivial linear ordinary differential
equations in all of the variables. For details see \cite{Draisma02}.
There are variants of this theorem that yield {\em polynomial}
realisations, such as the following.

\begin{thm}[Draisma, \cite{Draisma02}] \label{thm:Nilpotent}
Let $(\liea{g},\liea{h})$ be an effective transitive pair of
codimension $n$, and suppose that $\liea{g}$ has subalgebra $\liea{m}$
complementary to $\liea{h}$ such that for each element $X \in \liea{m}$
the linear map $\ad(X):\liea{g} \to \liea{g},\ Y \mapsto [X,Y]$ is
nilpotent. Then $(\liea{g},\liea{h})$ has a realisation in which all
coefficients are polynomial.
\end{thm}

In this case, the Realisation Formula with $Y_1,\ldots,Y_n$ any ordered
basis of $\liea{m}$ gives such a polynomial parameterisation.

\begin{ex}
In Morozov's Theorem~\ref{thm:Morozov} the second type of primitive
pairs are of the form described in this theorem. Indeed, assume that
$\liea{g}=\liea{h} \ltimes \liea{m}$ with $\liea{m}$ an Abelian ideal,
and choose a basis $Y_1,\ldots,Y_m$. Then we have $Y_iY_j=Y_jY_i$
in $U(\liea{g})$, so that the linear part of $Y_1^{\alpha_1} \cdots
Y_n^{\alpha_n} Y_i$ is only non-zero when all $\alpha_j$ are zero. This
shows that $Y_i$ is mapped to $\pafg{x_i}$ in the realisation formula.
On the other hand, for $X \in \liea{h}$ we have
\[ Y_j X=X Y_j + \sum_{i=1}^n a_{ij} Y_i \]
where $(a_{ij})_{ij}$ is the matrix of $-\ad_{\liea{m}} X$ with respect to
the basis $Y_1,\ldots,Y_n$. The linear part of a higher-degree monomial
in the $Y_i$ times $X$ is again zero. Hence we find that the realisation
maps $X$ to $\sum_j (\sum_i a_{ij} x_j \pafg{x_i})$.
\end{ex}

\begin{ex}
If $\liea{g}$ is a complex semisimple Lie algebra and $\liea{p}$ is a
parabolic subalgebra, then $\liea{p}$ has a complement as in this theorem,
and as a result the pair $(\liea{g},\liea{p})$ has a realisation with
polynomial coefficients. For explicit formulas, see \cite{Richter99}.
\end{ex}

\begin{re}
Table~\ref{tab:NonPrimitive} was compiled as follows: first, I re-did
the classification of effective, transitive pairs using only purely
Lie-algebraic arguments. It turns out that the primitive pairs are of the
form in Theorem~\ref{thm:Nilpotent}, while the non-primitive pairs are
trivially of the form in Theorem~\ref{thm:VeryImprimitive}. I then used
the pairs, with the appropriate basis of a complement, as input to a {\tt
GAP} implementation of the Realisation Formula
\cite{GAP4,blattner}, truncating the
power series at high enough degree to recognise the relevant polynomials
and exponentials (this, too, can be done automatically, because the
proofs of both theorems above gives explicit information on the degrees
of the polynomials and on which exponentials will appear).

As remarked earlier, Lie's {\em Gruppenregister} \cite{Lie24} does
not contain a complete list of transitive Lie algebras in three
variables. Indeed, the ones that are missing there are precisely the
very imprimitive ones, for which Theorem~\ref{thm:VeryImprimitive}
ensures the existence of a nice realisation. This fact, together with
Lie's classification of the other pairs, confirms Lie's quote above in
the case of three variables.  
\end{re}

There are pairs $(\liea{g},\liea{h})$ for which a realisation with
polynomial coefficients exists that can probably not be obtained directly
from the Realisation Formula. 

\begin{ex}
In Morozov's Theorem~\ref{thm:Morozov}, the first class of primitive
pairs $(\liea{k} \oplus \liea{k},\liea{h})$ (with $\liea{h}$ the diagonal
subalgebra) have polynomial realisations. A proof of this fact uses
the simply connected algebraic group $P$ with Lie algebra $\liea{k}$:
let $H$ be the diagonal subgroup of $P \times P$. Then the action of
$P \times P$ by left multiplication on the quotient $(P \times P)/H$
differentiates to a realisation of $\liea{p} \times \liea{p}$ by vector
fields on this quotient. More precisely, for every open affine subset
$U$ of the base point $(e,e)H$ (with stabiliser $H$) we obtain a Lie
algebra homomorphism $\liea{p} \times \liea{p} \to \Der K[U]$ such that
$\liea{h}$ is the preimage of all vector fields in $K[U]$ vanishing
at the base point. If we can choose $U$ to be isomorphic to an affine
{\em space} (and not just to some affine {\em variety}), then $K[U]$
is a polynomial ring and we have the desired polynomial realisation. By
the Bruhat decomposition, $P$ contains an open neighbourhood $V$ of $e$
that is isomorphic to an affine space, and we may take $U$ to be the image
of $V \times \{e\}$ in the quotient $P \times P / H$. This construction
generalises, in fact, to pairs corresponding to {\em spherical varieties}
\cite{Brion97}, \cite{Brion86}.
\end{ex}

Although I have no proof, I think that in the preceding
class of examples no choice of a basis complementary to $\liea{h}$
leads the Realisation Formula to output a polynomial realisation.
Observe that this is {\em not} in contradiction with the uniqueness
part in the Realisation Theorem: although any realisation of a pair
$(\liea{g},\liea{h})$ is related to a fixed realisation obtained from the
Realisation Formula by a unique formal coordinate change, this coordinate
change, which depends on countably many independent coefficients in the
formal power series images of $x_1,\ldots,x_n$, will typically {\em not}
come from a change of (basis of a) complement to $\liea{h}$, which is
determined by finitely many independent coefficients, and which is all
the data that goes into the Realisation Formula. This circumstance
and Lie's quote above lead to two interesting research problems:
\begin{itemize}

\item First, might Lie's observation literally hold true in higher
dimensions? In other words, does every transitive pair have a realisation
with coefficients that are polynomials in the variables $x_i$ and simple
exponentials $\exp(\lambda x_i)$? Might such a realisation always be
obtainable with the Realisation Formula?  Presently I am inclined to
answer both questions in the negative, but I have no counter-examples
yet. For instance, when allowing exponentials, the Realisation Formula
{\em does} yield allowed realisations for the ``diagonal pairs'' in the
previous example \cite{Draisma02}. My search for actual counter-examples
would start at primitive pairs $(\liea{g},\liea{h})$ from Dynkin's
list with $\liea{g}$ classical and its standard module an irreducible
$\liea{h}$-module.

\item Taking Lie's observation less literally, one might change the
condition 
\begin{enumerate}
\item that the coefficients of the vector fields be solutions to
linear o.d.e.s with constant coefficients in each of the variables
\end{enumerate}
to the condition 
\begin{enumerate}
\addtocounter{enumi}{1}
\item that the Lie algebra as a whole be the solution of a system of
linear p.d.e.s (not necessarily with constant coefficients).
\label{it:2}
\end{enumerate}
For example, the Lie algebra of type (8)
in Table~\ref{tab:Primitive} is the solution to the determining system
\eqref{eq:DetSys}, but in condition \eqref{it:2} we do not insist
that the system have constant coefficients.

This condition \eqref{it:2} was suggested by Mohamed Barakat at the 2010
Oberwolfach mini-workshop {\em Algebraic and analytic techniques for
polynomial vector fields}. There is a beautiful challenge here: start
with a transitive pairs $(\liea{g},\liea{h})$ and construct, preferably
in the coordinate-free manner introduced by Blattner, a system of linear
p.d.e.s to which (a realisation of) the pair is the solution. Then find
necessary and sufficient conditions for the system to have ``constant
coefficients''---which would correspond to the class described by Lie.
\end{itemize}

Especially the second research question seems very promising to me, and
with it I conclude this overview paper on Lie algebras of vector fields.

\end{document}